\documentclass[12pt,leqno]{article}
\usepackage{amsfonts,amsthm,amsmath}
\newtheorem{thm}{Theorem}[section]
\newtheorem{prop}[thm]{Proposition}
\newtheorem{lem}[thm]{Lemma}
\newtheorem{cor}[thm]{Corollary}
\theoremstyle{remark}
\newtheorem{rem}[thm]{Remark}
\newcommand{\FF}{\mathbb{F}}
\newcommand{\ZZ}{\mathbb{Z}}
\newcommand{\RR}{\mathbb{R}}

\newcommand{\allone}{\mathbf{1}}
\DeclareMathOperator{\supp}{supp}

\begin{document}

\title{Construction of 
Some Unimodular Lattices with Long Shadows}

\author{Masaaki Harada\thanks{Department of Mathematical Sciences,
Yamagata University, Yamagata 990--8560, Japan, and
PRESTO, Japan Science and Technology Agency (JST), Kawaguchi,
Saitama 332--0012, Japan. email: mharada@sci.kj.yamagata-u.ac.jp}}

\maketitle

\begin{center}
{\bf Dedicated to Professor Yasumasa Nishiura on His 60th Birthday}
\end{center}

\begin{abstract}
In this paper, we construct odd unimodular lattices in 
dimensions $n=36,37$ having minimum norm $3$ and 
$4 s=n-16$, where $s$ is the minimum norm of 
the shadow.
We also construct odd unimodular lattices in 
dimensions $n=41,43,44$ having minimum norm $4$ and 
$4 s=n-24$.
\end{abstract}

\section{Introduction}
Shadows of odd unimodular lattices appeared in~\cite{CS98} (see also 
\cite[p.~440]{SPLAG}),
and shadows play an important role in the study of odd unimodular lattices.
Let $L$ be an odd unimodular lattice in dimension $n$.
The {\em shadow} $S(L)$ of $L$ is defined to be $S(L)= L_0^* \setminus L$,
where $L_0$ denotes the even sublattice of $L$ and
$L_0^*$ denotes the dual lattice of $L_0$.
Note that the norm of a vector of $S(L)$ is congruent to $n/4$ modulo 
$2$~\cite{CS98}.

We define
\[
\sigma(L)=4 \min(S(L)),
\]
where $\min(S(L))$ denotes the minimum norm of $S(L)$.
Elkies~\cite{E95a} began the investigation of
odd unimodular lattices $L$ with long shadows, that is,
large $\sigma(L)$.
It was shown in~\cite{E95a} that $\sigma(L) \le n$
and $\ZZ^n$ is the only odd unimodular lattice $L$ with
$\sigma(L)=n$, up to isomorphism.
For the case $\sigma(L)<n$,
we may assume that $L$ has minimum norm at least $2$ 
(see~\cite{E95b,NV03}).
Elkies~\cite{E95b} determined all such lattices
for the case $\sigma(L)=n-8$.

For the next case $\sigma(L)=n-16$,
Nebe and Venkov~\cite{NV03} showed that if
there is an odd unimodular lattice $L$ with minimum
norm $\min(L) \ge 3$
then $n \le 46$.
Moreover, 
in this case, it is not known whether there is such a lattice with minimum
norm $3$
or not for only $n=36,37,38,39,41,42,43$
(see \cite[p.~148]{Ga07} and Table \ref{Tab:1}).
It follows from \cite[Table 1]{NV03}
that there is no odd unimodular lattice $L$ with $\min(L) \ge 4$
for the case $\sigma(L)=n-16$ (see also \cite[Table 2]{Gau07}).

For the case $\sigma(L)=n-24$, 
Gaborit~\cite{Ga07} showed that if
there is an odd unimodular lattice $L$ with $\min(L) = 4$
then $n \le 47$ (see also \cite[Table 2]{Gau07}).
Moreover, 
in this case, it is not known whether there is such a lattice
or not for only $n=37,41,43,44,45$
(see \cite[p.~148]{Ga07} and Table \ref{Tab:2}).

The aim of this paper is to provide the existence of
some unimodular lattices with long shadows
whose existences were previously not known.
These lattices are constructed from self-dual $\ZZ_k$-codes
($k=4,5$).
The paper is organized as follows.
In Section~\ref{Sec:2}, we give definitions and some basic
results on unimodular lattices and self-dual $\ZZ_k$-codes.
In Section \ref{Sec:36}, we give conditions to construct
odd unimodular lattices $L$ in dimension $36$
with $\min(L) =3$ and $\sigma(L)=20$ using
self-dual $\ZZ_4$-codes by Construction A 
(Proposition~\ref{prop:36con}).
By finding a self-dual $\ZZ_4$-code satisfying these
conditions, 
the first example of such a lattice
is constructed.
This dimension is the smallest dimension
of an odd unimodular lattice $L$ 
with $\min(L) =3$ and $\sigma(L)=n-16$
whose existence was previously not known.
A new odd unimodular lattice $L$ in dimension $36$
with $\min(L)=4$ and $\sigma(L)=12$ is constructed
as a neighbor of the above lattice.
In Section \ref{Sec:Ot},
by considering self-dual $\ZZ_k$-codes ($k=4,5$),
we construct an odd unimodular lattice $L$ in dimension $n$
with 
\begin{align*}
(n,\min(L),\sigma(L))=&
(37,3,n-16),
(41,4,n-24),
(43,4,n-24), \text{ and}
\\ &
(44,4,n-24).
\end{align*}
The current state of knowledge about the existence of
unimodular lattices $L$ in dimension $n$
with $(\min(L),\sigma(L))=
(3,n-16)$ and $(4,n-24)$ is listed 
in Tables \ref{Tab:1} and \ref{Tab:2}, respectively.
All computer calculations in this paper
were done by {\sc Magma}~\cite{Magma}.

 \section{Preliminaries}\label{Sec:2}

\subsection{Unimodular lattices}
A (Euclidean) lattice $L \subset \RR^n$ 
in dimension $n$
is {\em unimodular} if
$L = L^{*}$, where
the dual lattice $L^{*}$ of $L$ is defined as
$\{ x \in {\RR}^n \mid (x,y) \in \ZZ \text{ for all }
y \in L\}$ under the standard inner product $(x,y)$.
A unimodular lattice $L$ is {\em even} 
if the norm $(x,x)$ of every vector $x$ of $L$ is even
and {\em odd} otherwise.
An even unimodular lattice in dimension $n$
exists if and only if $n \equiv 0 \pmod 8$, while
an odd  unimodular lattice exists for every dimension.
The minimum norm $\min(L)$ of $L$ is the smallest 
norm among all nonzero vectors of $L$.
Two lattices $L$ and $L'$ are {\em isomorphic},
if there exists an orthogonal matrix $A$ with
$L' = L \cdot A =\{xA \mid x \in L\}$.
The automorphism group of $L$ is the group of all
orthogonal matrices $A$ with $L = L \cdot A$.

The theta series $\theta_{L}(q)$ of $L$ is the formal power
series
$
\theta_{L}(q) = \sum_{x \in L} q^{(x,x)}.
$
The kissing number is the second nonzero coefficient of the
theta series, that is, the number of vectors of minimum norm
in $L$.
Conway and Sloane~\cite{CS98} showed that
when the theta series of an odd unimodular lattice $L$
in dimension $n$
is written as
\begin{equation}
\label{eq:theta}
\theta_L(q)=
 \sum_{j =0}^{\lfloor n/8\rfloor} a_j\theta_3(q)^{n-8j}\Delta_8(q)^j,
\end{equation}
the theta series of the shadow $S(L)$
is written as
\begin{equation}
\label{eq:theta-S}
\theta_{S(L)}(q)= \sum_{j=0}^{\lfloor n/8\rfloor}
\frac{(-1)^j}{16^j} a_j\theta_2(q)^{n-8j}\theta_4(q^2)^{8j},
\end{equation}
where
$\Delta_8(q) = q \prod_{m=1}^{\infty} (1 - q^{2m-1})^8(1-q^{4m})^8$
and $\theta_2(q), \theta_3(q)$ and $\theta_4(q)$ are the Jacobi
theta series~\cite{SPLAG}.

\subsection{Self-dual $\ZZ_k$-codes and Construction A}

Let $\ZZ_{k}$ be the ring of integers modulo $k$, where $k$ 
is a positive integer greater than $1$.
A {\em $\ZZ_{k}$-code} $C$ of length $n$ 
is a $\ZZ_{k}$-submodule of $\ZZ_{k}^n$.
We shall exclusively deal with the 
case $k=4$.
Two $\ZZ_k$-codes are {\em equivalent} if one can be obtained from the
other by permuting the coordinates and (if necessary) changing
the signs of certain coordinates.
The {\em dual} code $C^\perp$ of $C$ is defined as
$C^\perp = \{ x \in \ZZ_{k}^n\ | \ x \cdot y = 0$ for all $y \in C\}$,
where $x \cdot y = x_1 y_1 + \cdots + x_n y_n$ for
$x=(x_1,\ldots,x_n)$ and $y=(y_1,\ldots,y_n)$.
A code $C$ is {\em self-dual} if $C=C^\perp.$
Let $C$ be a self-dual $\ZZ_k$-code of length $n$.
Then the following lattice
\[
A_{k}(C) = \frac{1}{\sqrt{k}}
\{(x_1,\ldots,x_n) \in \ZZ^n \mid
(x_1 \bmod k,\ldots,x_n \bmod k)\in C\}
\]
is a unimodular lattice in dimension $n$.
This construction of lattices is called Construction A.

\subsection{Self-dual $\ZZ_4$-codes}
\label{Subsec:Z4}
Let $C$ denote a $\ZZ_4$-code of length $n$.
The {\em Euclidean weight} of a codeword $x=(x_1,\ldots,x_n)$ of $C$ is
$n_1(x)+4n_2(x)+n_3(x)$, where $n_{\alpha}(x)$ denotes
the number of components $i$ with $x_i=\alpha$ $(\alpha=1,2,3)$.
The {\em minimum Euclidean weight} $d_E(C)$ of $C$ is the smallest Euclidean
weight among all nonzero codewords of $C$.
Every $\ZZ_4$-code $C$ of length $n$ has two binary codes 
$C^{(1)}$ and $C^{(2)}$ associated with $C$:
\[
C^{(1)}= \{ c \bmod 2 \mid  c \in C \} \text{  and }
C^{(2)}= \left\{ c \bmod 2 \mid c \in \ZZ_4^n, 2c\in C \right\}.
\]
The binary codes $C^{(1)}$ and $C^{(2)}$ are called the 
{\em residue} and {\em torsion} codes of $C$, respectively.
If $C$ is self-dual, then $ C^{(1)}$ is a binary 
doubly even code with $C^{(2)} = {C^{(1)}}^{\perp}$~\cite{Z4-CS}.
It is easy to see that 
$\min\{d(C^{(1)}),4d(C^{(2)}))\} \le d_E(C)$,
where $d(C^{(i)})$ denotes the minimum weight of $C^{(i)}$
$(i=1,2)$.
In addition, $d_E(C) \le 4d(C^{(2)})$ (see~\cite{H10}).
Also, 
$A_4(C)$ has minimum norm $\min\{4,d_E(C)/4\}$ (see~\cite{Z4-BSBM}).
Therefore, 
if $A_4(C)$ has minimum norm $3$ (resp.\ $4$),
then $C$ must have minimum Euclidean weight $12$ 
(resp.\ at least $16$) and
$C^{(2)}$ has minimum weight at least $3$ (resp.\ at least $4$).


Two $\ZZ_4$-codes 
differing by only a permutation of coordinates are called
permutation-equivalent.
Any self-dual $\ZZ_4$-code $C$ of length $n$
with residue code of dimension $k_1$
is permutation-equivalent to a code
$C'$ with generator matrix
in standard form
\begin{equation}
\label{eq:g-matrix}
\left(\begin{array}{ccc}
I_{k_1} & A & B_1+2B_2 \\
O    &2I_{n-2k_1} & 2D 
\end{array}\right),
\end{equation}
where $A$, $B_1$, $B_2$, $D$ are $(1,0)$-matrices,
$I_k$ denotes the identity matrix of order $k$, and
$O$ denotes the zero matrix~\cite{Z4-CS}.
In this paper, 
when we give a generator matrix of $C$, 
we consider a generator matrix
in standard form (\ref{eq:g-matrix})
of $C'$, which is permutation-equivalent to $C$,
then we only list the 
$k_1 \times (n-k_1)$ matrix
$
\left(\begin{array}{cc}
 A & B_1+2B_2 
\end{array}\right)
$
to save space.
Note that 
$\left(\begin{array}{ccc}
O    &2I_{n-2k_1} & 2D 
\end{array}\right)$
in (\ref{eq:g-matrix})
can be obtained from 
$\left(\begin{array}{ccc}
I_{k_1} & A & B_1+2B_2 
\end{array}\right)$
since $C'^{(2)} = {C'^{(1)}}^{\perp}$.

\section{Dimension $n=36$ and $\sigma(L)=n-16, n-24$}\label{Sec:36}
In this section, we give conditions to construct
odd unimodular lattices $L$ in dimension $n=36$
with $\min(L) =3$ and $\sigma(L)=n-16$ from
self-dual $\ZZ_4$-codes by Construction A.
By finding a self-dual $\ZZ_4$-code satisfying these
conditions, the first example of such a lattice
is constructed.
A new optimal odd unimodular lattice $L$ 
with $\sigma(L)=n-24$ is also constructed from the above lattice.

Let $L_{36}$ be an odd unimodular lattice 
in dimension $36$ having minimum norm at least $3$.
Using (\ref{eq:theta}) and (\ref{eq:theta-S}),
it is easy to determine
the possible theta series $\theta_{L_{36}}(q)$ and
$\theta_{S(L_{36})}(q)$ of $L_{36}$ and its shadow $S(L_{36})$:
\begin{align*}
\theta_{L_{36}}(q) =&
1 
+ (960 - \alpha)q^3 
+ (42840 + 4096 \beta)q^4 
\\ & \hspace{4cm}
+ (1882368 + 36 \alpha - 98304 \beta)q^5
+ \cdots,
\\
\theta_{S(L_{36})}(q) =&
\beta q  + (\alpha  - 60 \beta) q^3 
+ (3833856 - 36 \alpha + 1734 \beta)q^5
+ \cdots,
\end{align*}
respectively,
where $\alpha$ and $\beta$ are nonnegative integers.
Then, the following lemma is immediate.

\begin{lem}\label{lem:36}
Let $L$ be an odd unimodular lattice in dimension $36$
having minimum norm $3$. 
Then, the kissing number of $L$ is at most $960$, 
and the equality holds if and only if
$\sigma(L)=20$. 
\end{lem}

Now, we give a method for construction of 
unimodular lattices $A_4(C)$ with $\min(A_4(C))=3$ and 
$\sigma(A_4(C))=20$,
using self-dual $\ZZ_4$-codes $C$.
It is known that a binary $[36,k,3]$ code exists only if $k \le 30$
(see~\cite{Brouwer-Handbook}). 
Hence, the dimension of the residue code of 
a self-dual $\ZZ_4$-code of length $36$
and minimum Euclidean weight $12$ is at least $6$.

\begin{prop}\label{prop:36con}
Let $C$ be a self-dual $\ZZ_4$-code of length $36$ such that
$C^{(1)}$ is a binary doubly even $[36,6,16]$ code and 
$C^{(2)}$ has minimum weight $3$.
Then $A_4(C)$ is a unimodular lattice with $\min(A_4(C))=3$
and $\sigma(A_4(C))=20$. 
\end{prop}
\begin{proof}
Let $B$ be a binary doubly even $[36,6,16]$ code
such that $B^\perp$ has minimum weight $3$.
The weight enumerator of $B$ is written as:
\[
W_{B}(y)=
1
+a y^{16}
+b y^{20}
+c y^{24}
+d y^{28}
+e y^{32}
+(2^6-1-a-b-c-d-e) y^{36},
\]
where $a,b,c,d$ and $e$ are nonnegative integers.
Since $B^\perp$ contains a codeword of weight $3$,
$2^6-1-a-b-c-d-e=0$.
By the MacWilliams identity, 
the weight enumerator of $B^\perp$
is given by:
\begin{align*}
W_{B^\perp}(y)=
&1 
+ \frac{1}{8}
(- 216 + 4a + 3b + 2c + d ) y 
+ (378 -6a - 6b - 5c - 3d ) y^2 
\\&
+ \frac{1}{8}
(- 24024 + 388a +403b + 386c+ 273d ) y^3 +
\cdots.
\end{align*}
Since $B^\perp$ has minimum weight $3$, we have 
\[
c= 270 - 6 a - 3 b \text{ and }
d=-324 + 8 a + 3 b.
\]
Thus, the weight enumerators of $B$ and $B^\perp$
are written using $a$ and $b$:
\begin{align*}
W_{B}(y)=&
1 
+ a y^{16} 
+ b y^{20} 
+ (270 -6a - 3b)y^{24} 
+ (-324 + 8a + 3b) y^{28}
\\ &
+ (117 -3a - b) y^{32},
\\
W_{B^\perp}(y)=&
1 + (- 1032 + 32a + 8b ) y^3 
+ ( 17649 -448 a - 128 b) y^4 + \cdots,
\end{align*}
respectively.
Then, only $(a,b)=(27,36)$ satisfies the condition that
the above coefficients in $W_{B}(y)$ and $W_{B^\perp}(y)$ 
are nonnegative integers.
Hence, the weight enumerators of $B$ and $B^\perp$ are uniquely 
determined as 
\begin{align}
\label{eq:WE36}
W_{B}(y)=&
1+  27 y^{16} +  36 y^{20},
\\
\label{eq:WE36d}
W_{B^\perp}(y)=&
1 + 120 y^3 + 945y^4 + 5832 y^5 +  30576 y^6 + 130680 y^7 + \cdots,
\end{align}
respectively.


As described in Section \ref{Subsec:Z4}, we have
\begin{equation}\label{eq:dE}
\min\{d(C^{(1)}),4d(C^{(2)}))\} \le d_E(C) \le 4d(C^{(2)}). 
\end{equation}
Hence, $d_E(C)=12$ and $A_4(C)$ has minimum norm $3$.
Let $e_i$ denote the $i$-th unit vector
$(\delta_{i,1},\delta_{i,2},\ldots,\delta_{i,36})$,
for $i=1,2,\ldots,36$, where $\delta_{ij}$ is the Kronecker delta.
By (\ref{eq:WE36d}), there are $120$ codewords of weight $3$
in $C^{(2)}$, and 
we denote the set of the $120$ codewords
by $C^{(2)}_3$.
Then, $A_4(C)$ contains the following set of vectors of norm $3$:
\[
\{\pm e_{j_1} \pm e_{j_2} \pm e_{j_3} \mid 
\{j_1,j_2,j_3\} \in S\},
\]
where $S=\{\supp(x) \mid x \in C^{(2)}_3\}$,
and $\supp(x)$ denotes the support of $x$.
Hence, there are at least $960$ vectors of norm $3$ in $A_4(C)$.
By Lemma~\ref{lem:36}, 
the result follows.
\end{proof}

It was shown in~\cite{PST} that there are four inequivalent
binary $[36,7,16]$ codes containing the all-one vector $\allone$.
Such a code and its dual code have
the following weight enumerators:
\begin{align*}
&1+ 63 y^{16}+ 63 y^{20}+ y^{36},
\\
&
1
+    945 y^{4}
+   30576 y^{6}
+  471420 y^{8}
+ 3977568 y^{10}
+ \cdots,
\end{align*}
respectively.
Let $B_{36,7,1}$ denote the
binary $[36,7,16]$ code containing $\allone$ with
automorphism group of order $1451520$
which is the symplectic group $Sp(6,2)$.
Using an approach used in \cite[Section~4]{HLM},
we verified that $B_{36,7,1}$ contains only one
doubly even $[36,6,16]$ subcode $B_{36,6}$ such that
$B_{36,6}^\perp$ has minimum weight $3$ up to equivalence.
The weight enumerators of $B_{36,6}$ and $B_{36,6}^\perp$
are given by (\ref{eq:WE36}) and (\ref{eq:WE36d}), respectively.
We verified  by {\sc Magma} 
that $B_{36,6}$ has automorphism group of  
order $51840$ containing
the symplectic group $Sp(4,3)$ as a subgroup of index $2$.

Starting from a given binary doubly even code $B$,
a method for construction of all
self-dual $\ZZ_4$-codes $C$ with $C^{(1)}=B$
was given in~\cite[Section 3]{Z4-PLF}.
Using this method,
we construct a self-dual $\ZZ_4$-code $C_{36}$ 
with $C_{36}^{(1)}=B_{36,6}$ explicitly.
By Proposition~\ref{prop:36con},
$A_4(C_{36})$ is the desired
unimodular lattice with $\sigma(A_4(C_{36}))=20$,
and we have the following:

\begin{prop}
There is a unimodular lattice $L$ in dimension $36$ having
minimum norm $3$ with $\sigma(L)=20$.
\end{prop}

We verified by {\sc Magma} that the unimodular lattice
$A_4(C_{36})$ has automorphism group of
order $1698693120$.
For the code $C_{36}$, we give a generator matrix
in standard form (\ref{eq:g-matrix}), by only listing 
the $6 \times 30$ matrix:
\[
\left(\begin{array}{cc}
 A & B_1+2B_2 
\end{array}\right)
=
\left(\begin{array}{cc}
 000010001011100111100010 &130113\\
 011001000101000011101001 &231133\\
 100010010111110001110111 &000230\\
 110011100011001100111111 &032133\\
 011110101100111001011001 &321113\\
 111111111111110000000000 &222203
\end{array}\right).
\]

There are three unimodular lattices 
containing the even sublattice of $A_4(C_{36})$.
We denote the two unimodular lattices rather than
$A_4(C_{36})$ by $N_{36}$ and $N'_{36}$.
It follows from the theta series of 
$A_4(C_{36})$ and $S(A_4(C_{36}))$
that both $N_{36}$ and $N'_{36}$ have 
minimum norm $4$ and kissing number $42840$,
thus, 
$\sigma(N_{36})=\sigma(N'_{36})=16$.
We verified by {\sc Magma} that the two lattices 
are isomorphic.
We also verified by {\sc Magma} that $N_{36}$ has
automorphism group of order $849346560$,
which is different to those of
two previously known unimodular lattices with
minimum norm $4$ and kissing number $42840$
in~\cite{lattice-datebase}.

Let $C_{36,0}$ be the subcode of $C_{36}$ consisting 
of codewords of Euclidean weight divisible by $8$.
Then $C_{36,0}$ is a subcode of index $2$ in $C_{36}$
(see \cite[Lemma 3.1]{DHS01}).
By Proposition 3.8 in~\cite{DHS01}, 
there is a self-dual $\ZZ_4$-code $D_{36}$ 
containing $C_{36,0}$ with $A_4(D_{36}) = N_{36}$.
For the code $D_{36}$, 
we give a generator matrix 
in standard form (\ref{eq:g-matrix}), by only listing 
the $6 \times 30$ matrix:
\[
\left(\begin{array}{cc}
 A & B_1+2B_2 
\end{array}\right)
=
\left(\begin{array}{cc}
 0001101111001000100110& 0131310\\
 1101111000101001010010& 1233122\\
 0111100011110010100011& 0121203\\
 1110100110001110111011& 3133300\\
 0101001101100001101101& 1203301\\
 0000011111111111111110& 0022113\\
 1111111111111100000000& 2220232
\end{array}\right).
\]
We verified that the residue code $D_{36}^{(1)}$ 
is equivalent to $B_{36,7,1}$.

\begin{rem}
The remaining three
binary $[36,7,16]$ codes containing $\allone$
have automorphism groups of orders
$10752$, $1920$ and $672$~\cite{PST}.
We denote these codes by $B_{36,7,2}$,
$B_{36,7,3}$ and $B_{36,7,4}$, respectively.
We verified that
$B_{36,7,2}$, $B_{36,7,3}$ and $B_{36,7,4}$
contain $1,2$ and $1$ 
doubly even $[36,6,16]$ subcodes such that
the dual codes have minimum weight $3$, respectively,
up to equivalence, and
the four codes and $B_{36,6}$ are inequivalent to each other.
By Proposition~\ref{prop:36con},
the four inequivalent codes 
rather than $B_{36,6}$ also 
give examples of unimodular lattices $L$ with 
$\min(L)=3$ and $\sigma(L)=20$.
\end{rem}

\section{Other cases}\label{Sec:Ot}
\subsection{Dimension $n=37$ and $\sigma(L)=n-16$}

Let $L_{37}$ be an odd unimodular lattice 
in dimension $n=37$ having 
minimum norm at least $3$.
We give the possible theta series 
of $L_{37}$ and its shadow $S(L_{37})$:
\begin{align*}
\theta_{L_{37}}(q) =&
1 
+ (1184 - \alpha)q^3 
+ (37962 - 2 \alpha + 2048 \beta) q^4 
\\ & \hspace{4cm}
+ (1758240 + 36 \alpha - 45056 \beta)q^5 
+ \cdots,
\\
\theta_{S(L_{37})}(q) =&
\beta q^{5/4} 
+ (2 \alpha  - 59 \beta)q^{13/4} 
+ (8486912 - 70 \alpha + 1674 \beta)q^{21/4}
+ \cdots,
\end{align*}
respectively,
where $\alpha$ and $\beta$ are nonnegative integers.
It turns out that $L_{37}$ has minimum norm $3$ and
kissing number $1184$ 
if and only if $\sigma(L_{37})=21$. 


It is known that a binary $[37,k,3]$ code exists only if $k \le 31$
(see~\cite{Brouwer-Handbook}). 
Hence, the dimension of the residue code of 
a self-dual $\ZZ_4$-code of length $37$
and minimum Euclidean weight $12$ is at least $6$.
When $n=37$, we have a weaker result than
Proposition~\ref{prop:36con}.

\begin{prop}
Let $B$ be a binary doubly even $[37,6]$
code such that the dual code has minimum weight $3$.
Then the weight enumerator $W_B(y)$ of $B$ is
given by:
\begin{align}
\label{eq:WE371}
W_B(y)        &= 1 + 20 y^{16} + 42 y^{20} + y^{24} \text{ or}, \\
\label{eq:WE372}
W_B(y)        &= 1+ y^{12} + 17 y^{16} + 45 y^{20}. 
\end{align}
\end{prop}
\begin{proof}
The proof is similar to that of Proposition~\ref{prop:36con}.
Let $B$ be a binary doubly even $[37,6]$ code
such that $B^\perp$ has minimum weight $3$.
The weight enumerator $W_B(y)$ is written as:
\begin{multline*}
W_{B}(y)=
1
+ a y^4 
+ b y^8 
+ c y^{12} 
+ d y^{16} 
+ e y^{20} 
+ f y^{24} 
+ g y^{28} 
\\
+ h y^{32} 
+ (2^6-1-a - b - c - d - e - f - g - h) y^{36},
\end{multline*}
where $a,b,c,d,e,f,g$ and $h$ are nonnegative integers.
Then the weight enumerator of $B^\perp$
is given by:
\begin{align*}
W_{B^\perp}(y)=&
1 
+ \frac{1}{8}
(- 271 + 8a + 7 b + 6 c + 5 d + 4 e + 3 f + 2 g +  h )y 
\\&
+ \frac{1}{8}
(4761 -24 a - 49 b - 66 c - 75 d - 76 e - 69 f - 54 g - 31 h )y^2 
\\&
+ \frac{1}{8}
(- 50295 + 1256 a + 959 b + 830 c + 805 d + 820 e + 811 f  
\\ &
\qquad 
+ 714 g + 465 h )y^3 + \cdots.
\end{align*}
From the condition that $B^\perp$ has minimum weight $3$,
we have 
\begin{align*}
g&= 455 - 28a - 21b - 15c - 10d - 6e - 3f,\\
h&= -639 + 48a + 35b + 24c + 15d + 8e + 3f.
\end{align*}
Thus, the weight enumerators are written as:
\begin{align*}
W_{B}(y)=&
1 + a y^4 + b y^8 + c y^{12} + d y^{16} + e y^{20} 
+ fy^{24} 
\\ &+  (455- 28a - 21b - 15c - 10d - 6e - 3f) y^{28} 
\\ &+  (-639 + 48a + 35b + 24c + 15d + 8e + 3f ) y^{32} 
\\ & + (247 -21a - 15b -10c - 6d - 3e - f )y^{36},
\\
W_{B^\perp}(y)=&
1 + (- 2820 + 448a + 280b + 160c + 80d + 32e + 8f )y^3 
+ \cdots.
\end{align*}
Then, only
\[
(a,b,c,d,e,f)=
(0, 0, 0, 20, 42, 1)\text{ and } (0, 0, 1, 17, 45, 0)
\]
satisfy the condition that the coefficients in 
$W_B(y)$ are nonnegative integers.
These cases correspond to
(\ref{eq:WE371}) and (\ref{eq:WE372}), respectively.
\end{proof}

For (\ref{eq:WE371}) and (\ref{eq:WE372}),
$W_{B^\perp}(y)$ is given by:
\begin{align}
\label{eq:WE371d}
W_{B^\perp}(y)&= 1 + 132 y^3 + 1072 y^4 + 6705 y^5 + 36324 y^6 + \cdots, \\
\label{eq:WE372d}
W_{B^\perp}(y)&= 1+ 140 y^3 + 1080 y^4 + 6633 y^5 + 36252 y^6 + \cdots,
\end{align}
respectively.

\begin{cor}\label{cor:37con}
Let $C$ be a self-dual $\ZZ_4$-code of length $37$  such that
$C^{(1)}$ is a binary doubly even $[37,6,12]$ code and
$C^{(2)}$ has minimum weight $3$.
Then $A_4(C)$ is a unimodular lattice with $\min(A_4(C))=3$ 
and kissing number at least $1120$.
\end{cor}
\begin{proof}
By (\ref{eq:dE}), $C$ has minimum Euclidean weight $12$.
It follows from (\ref{eq:WE372d}) that
$A_4(C)$ has at least $1120$ vectors of norm $3$.
\end{proof}


As a subcode of some binary maximal doubly even code of length $37$,
we found a doubly even $[37,6,12]$ code $B_{37}$ such that
$B_{37}^\perp$ has minimum weight $3$.
The weight enumerators of $B_{37}$ and $B_{37}^\perp$
are given by (\ref{eq:WE372}) and (\ref{eq:WE372d}), respectively.
We verified  by {\sc Magma} 
that $B_{37}$ has automorphism group of  
order $120$ which is the symmetric group of degree $5$.

Corollary~\ref{cor:37con} only guarantees that
$A_4(C)$ has minimum norm $3$ and
kissing number at least 
$1120$ for 
a self-dual $\ZZ_4$-code $C$ with $C^{(1)}=B_{37}$.
Using the method in~\cite[Section 3]{Z4-PLF},
we found a self-dual $\ZZ_4$-code $C_{37}$ 
such that $C_{37}^{(1)}=B_{37}$ and 
the kissing number of $A_4(C_{37})$ is exactly $1184$.
Therefore, we have the following:

\begin{prop}
There is a unimodular lattice $L$ in dimension $37$ having
minimum norm $3$ with $\sigma(L)=21$.
\end{prop}

We verified by {\sc Magma} that the unimodular lattice
$A_4(C_{37})$ has automorphism group of
order $7864320$.
For the code $C_{37}$, we give a generator matrix
in standard form (\ref{eq:g-matrix}), by only listing 
the $6 \times 31$ matrix:
\[
\left(\begin{array}{cc}
 A & B_1+2B_2 
\end{array}\right)
=
\left(\begin{array}{cc}
 0101001010110011110000011 & 003121\\
 0111100111010100000110110 & 323200\\
 1010001100010011111110101 & 313301\\
 0000001110001100111101100 & 033031\\
 1110011111111010100000000 & 000132\\
 1001110010000111111111000 & 220010
\end{array}\right).
\]

\subsection{Dimension $n=41$ and $\sigma(L)=n-24$}

Let $L_{41}$ be an odd unimodular lattice 
in dimension $41$ having minimum norm $4$.
We give the possible theta series 
of $L_{41}$ and its shadow $S(L_{41})$:
\begin{align*}
\theta_{L_{41}}(q) =&
1 
+ (15170 + 128\alpha )q^4
+ (1226720 - 1792\alpha - 524288\beta) q^5
+ \cdots,
\\
\theta_{S(L_{41})}(q) =&
\beta q^{1/4}
+ (\alpha - 79\beta) q^{9/4}
+ (104960 - 55\alpha + 3040\beta)q^{17/4}
+ \cdots,
\end{align*}
respectively,
where $\alpha$ and $\beta$ are nonnegative integers.
It turns out that $L_{41}$ has kissing number $15170$ 
if and only if $\sigma(L_{41})=17$. 


Let $C_{41}$ be the $\ZZ_4$-code with generator matrix
in standard form (\ref{eq:g-matrix}), where
$\left(\begin{array}{cc}
 A & B_1+2B_2 
\end{array}\right)$ is listed in Figure \ref{Fig:41}.
We verified that $C_{41}$ is a self-dual $\ZZ_4$-code of
minimum Euclidean weight $16$ such that
$A_4(C_{41})$ has kissing number $15170$.
Hence, we have the following:

\begin{prop}
There is a unimodular lattice $L$ in dimension $41$ having
minimum norm $4$ with $\sigma(L)=17$.
\end{prop}

\begin{figure}[thb]
\centering
{\small
\[
\left(\begin{array}{cc}
 A & B_1+2B_2 
\end{array}\right)
=
\left(\begin{array}{cc}
 10111111111001000010001& 101031013\\
 11001100111111100111000& 112233010\\
 01101100001100110110100& 001332012\\
 11010101000111001010100& 023102310\\
 11111010110101011110111& 033211301\\
 10000101001101101010101& 132003201\\
 11100110011010001101111& 210211330\\
 01110001111111010000111& 230321312\\
 01101110010011001111011& 232022022
\end{array}\right)
\]
\caption{Generator matrix of $C_{41}$}
\label{Fig:41}
}
\end{figure}

The residue code $C_{41}^{(1)}$ is a binary doubly even
$[41,9,12]$ code with a trivial automorphism group
and weight enumerator:
\[
1
+     y^{12}
+  89 y^{16}
+ 288 y^{20}
+ 108 y^{24}
+  23 y^{28}
+   2 y^{32}.
\]

\subsection{Dimension $n=43$ and $\sigma(L)=n-24$}

Let $L_{43}$ be an odd unimodular lattice 
in dimension $43$ having minimum norm $4$.
We give the possible theta series 
of $L_{43}$ and its shadow $S(L_{43})$:
\begin{align*}
\theta_{L_{43}}(q) =&
1 
+ (9030 + 32 \alpha)q^4 
+ (941184 - 320\alpha - 131072\beta)q^5 
+ \cdots,
\\
\theta_{S(L_{43})}(q) =&
\beta q^{3/4} 
+ (\alpha - 77\beta)q^{11/4} 
+ (660480 - 53 \alpha + 2883 \beta)q^{19/4} 
+ \cdots,
\end{align*}
respectively,
where $\alpha$ and $\beta$ are nonnegative integers.
It turns out that $L_{43}$ has kissing number $9030$ 
if and only if $\sigma(L_{43})=19$. 

Let $C_{43}$ be the $\ZZ_4$-code with generator matrix
in standard form (\ref{eq:g-matrix}), 
where 
$\left(\begin{array}{cc}
 A & B_1+2B_2 
\end{array}\right)$ is listed in Figure \ref{Fig:43}.
We verified that $C_{43}$ is a self-dual $\ZZ_4$-code of
minimum Euclidean weight $16$ such that
$A_4(C_{43})$ has kissing number $9030$.
Hence, we have the following:

\begin{prop}
There is a unimodular lattice $L$ in dimension $43$ having
minimum norm $4$ with $\sigma(L)=19$.
\end{prop}

\begin{figure}[thb]
\centering
{\small
\[
\left(\begin{array}{cc}
 A & B_1+2B_2 
\end{array}\right)
=
\left(\begin{array}{cc}
01010011110100111 &0033202030013\\
01100001111010111 &3212001202013\\
10100011110100010 &1122313032032\\
00001111111011100 &1100211211331\\
01101101000000110 &3300133330012\\
00011011001111100 &2213322021123\\
00111011110010110 &1230201103220\\
00100000101001101 &0121022213003\\
11000010110100011 &1012310010101\\
01011101010110101 &3321023022100\\
10001011010111000 &3021302111320\\
00001000111100101 &3231010323110\\
11111111111111000 &0222002002232
\end{array}\right)
\]
\caption{Generator matrix of $C_{43}$}
\label{Fig:43}
}
\end{figure}

The residue code $C_{43}^{(1)}$ is a binary doubly even
$[43,13,12]$ code with a trivial automorphism group
and weight enumerator:
\[
1
+   29 y^{12}
+ 1067 y^{16}
+ 3498 y^{20}
+ 3010 y^{24}
+  569 y^{28}
+   18 y^{32}.
\]

\subsection{Dimension $n=44$ and $\sigma(L)=n-24$}

Let $L_{44}$ be an odd unimodular lattice 
in dimension $44$ having minimum norm $4$.
We give the possible theta series 
of $L_{44}$ and its shadow $S(L_{44})$:
\begin{align*}
\theta_{L_{44}}(q) =&
1 
+ (6600 + 16 \alpha) q^4 
+ (811008  - 128 \alpha- 65536 \beta) q^5 
+ \cdots,
\\
\theta_{S(L_{44})}(q) =&
\beta q
+(\alpha - 76 \beta) q^3
+ (1622016 - 52 \alpha + 2806 \beta) q^5 
+ \cdots,
\end{align*}
respectively,
where $\alpha$ and $\beta$ are nonnegative integers.
It turns out that $L_{44}$ has kissing number $6600$ 
if and only if $\sigma(L_{44})=20$. 

Let $C_{44}$ be the $\FF_5$-code of length $44$
whose generator matrix is
\[
\left(
\begin{array}{ccc@{}c}
\quad & {\Large I_{22}} & \quad &
\begin{array}{cc}
A & B \\
-B^T & A^T
\end{array}
\end{array}
\right),
\]
where 
$A$ and $B$ are $11 \times 11$ 
negacirculant matrices with the first rows
\[
(1, 0, 0, 0, 0, 3, 3, 0, 1, 0, 3) \text{ and }
(0, 1, 1, 0, 2, 1, 0, 0, 0, 2, 3),
\]
respectively, 
that is, $A$ and $B$ have the following form
\[
\left( \begin{array}{ccccc}
r_1&r_2&r_3& \cdots &r_{11} \\
-r_{11}&r_1&r_2& \cdots &r_{10} \\
-r_{10}&-r_{11}&r_1& \cdots &r_{9} \\
\vdots &\vdots & \vdots && \vdots\\
-r_2&-r_3&-r_4& \cdots&r_1
\end{array}
\right),
\]
and $A^T$ denotes the transposed matrix of $A$.
Since $AA^T+BB^T=4I_{11}$,
$C_{44}$ is self-dual.
In addition, it can be checked that
$C_{44}$ has the following weight enumerator:
\[
1+ 924 y^{12}+  1848 y^{13} 
+ 19272  y^{14}  +91168 y^{15} 
+\cdots.
\]
We verified by {\sc Magma} that
the unimodular lattice $A_5(C_{44})$ has
minimum norm $4$ and kissing number $6600$.
Hence, we have the following:

\begin{prop}
There is a unimodular lattice $L$ in dimension $44$ having
minimum norm $4$ with $\sigma(L)=20$.
\end{prop}

Our computer search failed to discover a 
unimodular lattice with long shadow
using a self-dual $\ZZ_4$-code for this case
and the remaining cases that
the existences are not known.

\subsection{Summary}
As a summary,  we list the number $\#$ of known non-isomorphic
unimodular lattices $L$ in dimension $n$
with $\min(L)=3$ and $\sigma(L)=n-16$
(resp.\ $\min(L)=4$ and $\sigma(L)=n-24$)
in Table \ref{Tab:1} (resp.\ Table \ref{Tab:2}).
Both tables update the two tables given in \cite[p.~148]{Ga07}.
We remark that the existence of
a unimodular lattice $L$ in dimension $37$ having minimum norm $4$
is still unknown for any $\sigma(L)$ (see \cite{lattice-datebase}).

\begin{table}[th]

\centering
\caption{Unimodular lattices $L$ with $\min(L)=3$ and $\sigma(L)=n-16$}
\vspace*{0.2in}
\label{Tab:1}
{\small
\begin{tabular}{c|c|l||c|c|l}
\noalign{\hrule height1pt}
$n$ & $\#$ & \multicolumn{1}{c||}{References} &
$n$ & $\#$ & \multicolumn{1}{c}{References} \\
\hline
23& $1$     & (see \cite{NV03}) &35& $\ge 1$ & \cite{NV03}  \\
24& $1$     & (see \cite{NV03}) &36& $\ge 1$ & $A_4(C_{36})$ in Section~\ref{Sec:36} \\
25& $0$     & (see \cite{NV03}) &37& $\ge 1$ & $A_4(C_{37})$ in Section~\ref{Sec:Ot} \\
26& $1$     & (see \cite{NV03}) &38& ? & \\
27& $2$     & (see \cite{NV03}) &39& ? & \\
28& $36$    & (see \cite{NV03}) &40& $\ge 1$ & (see \cite{Ga07}) \\
29& $\ge 1$ & \cite{NV03}  &41& ? & \\
30& $\ge 1$ & \cite{NV03}  &42& ? & \\
31& $\ge 1$ & \cite{NV03}  &43& ? & \\
32& $\ge 1$ & \cite{NV03}  &44& 0 & \cite{NV03} \\
33& $\ge 1$ & \cite{NV03}  &45& 0 & \cite{NV03} \\
34& $\ge 1$ & \cite{NV03}  &46& 1 & \cite{NV03} \\
\noalign{\hrule height1pt}
\end{tabular}
}
\end{table}

\begin{table}[th]
\centering
\caption{Unimodular lattices $L$ with $\min(L)=4$ and $\sigma(L)=n-24$}
\vspace*{0.2in}
\label{Tab:2}
{\small
\begin{tabular}{c|c|l||c|c|l}
\noalign{\hrule height1pt}
$n$ & $\#$ & \multicolumn{1}{c||}{References} &
$n$ & $\#$ & \multicolumn{1}{c}{References} \\
\hline
32 & 5       & \cite{CS98} &42 & $\ge 1$ & (see \cite{Ga07}) \\
36 & $\ge 3$ & \cite{lattice-datebase}, $A_4(D_{36})$ in Section~\ref{Sec:36} 
        &43 & $\ge 1$ & $A_4(C_{43})$ in Section~\ref{Sec:Ot}  \\
37 & ?       &                   
        &44 & $\ge 1$ & $A_5(C_{44})$ in Section~\ref{Sec:Ot}      \\
38 & $\ge 1$ & (see \cite{Ga07}) &45 & ?       &                   \\
39 & $\ge 1$ & (see \cite{Ga07}) &46 & $\ge 1$ & (see \cite{Ga07}) \\
40 & $\ge 1$ & (see \cite{Ga07}) &47 & $\ge 1$ & (see \cite{Ga07}) \\
41 & $\ge 1$ & $A_4(C_{41})$ in Section~\ref{Sec:Ot} &&&\\
\noalign{\hrule height1pt}
\end{tabular}
}
\end{table}

\bigskip
\noindent
{\bf Acknowledgment.} 
The author would like to thank the anonymous referee
for helpful suggestions that improved
Proposition \ref{prop:36con} and 
Corollary \ref{cor:37con}.
This work was supported by JST PRESTO program.


\end{document}